# Rewriting Systems and Geometric 3-Manifolds

*Susan Hermiller and Michael Shapiro*

**Abstract:** The fundamental groups of most (conjecturally, all) closed 3-manifolds with uniform geometries have finite complete rewriting systems. The fundamental groups of a large class of amalgams of circle bundles also have finite complete rewriting systems. The general case remains open.

## §1. Introduction

In this paper we point out that well-known properties of finite complete rewriting systems and well-known facts about geometric 3-manifolds combine to give the following. (See below for definitions.)

**Theorem 1.** *Suppose that $M$ is a closed 3-manifold bearing one of Thurston's eight geometries. Suppose further that if $M$ is hyperbolic, that $M$ virtually fibers over a circle. Then $\pi_1(M)$ has a finite complete rewriting system.*

According to a conjecture of Thurston (**[Th]**, question 18), every closed hyperbolic 3-manifold obeys the last hypothesis.

We also exhibit a class of non-uniform geometric 3-manifolds whose fundamental groups have finite complete rewriting systems. In particular, suppose that $M$ is a graph of circle bundles based on a graph $\Gamma$. We will call an edge of this graph a ***loop*** if it has the same initial and terminal vertex. We suppose that when all loops are removed, the resulting graph is a tree. Under certain conditions on the way the vertex manifolds are glued along their boundary tori, the fundamental group $\pi_1(M)$ has a finite complete rewriting system.

AMS Classifications 20F32, 68Q42, 57M05.

The first author wishes to thank the National Science Foundation for partial support from grant DMS-923088.

The second author wishes to thank the Australian Research Council.

## §2. Proof of Theorem 1

We review the appropriate background and definitions.

Let $G$ be a group with finite generating set $A$. We write $A^*$ for the free monoid on $A$. Each element of $A$ evaluates into $G$ under the identity map and this extends to a unique monoid homomorphism of $A^*$ onto $G$ which we denote by $w \mapsto \overline{w}$.

A **rewriting system** $R$ over the set $A$ is a subset of $A^* \times A^*$. We write a pair $(u, v) \in R$ as $u \to v$ and call this a **rewriting rule** or *replacement rule*. If $u \to v$ is a rewriting rule, then for any $xuy \in A^*$, we write $xuy \to xvy$.

We say a finite set $R = \{u_i \to v_i\}$ is a **finite complete rewriting system** for $(G, A)$ if
1) The monoid presentation $\langle A \mid u_i = v_i \rangle$ is a presentation of the underlying monoid of G,
2) For each element $g \in G$ there is exactly one word $w \in A^*$ so that $g = \overline{w}$ and $w$ contains no $u_i$ as a substring (that is, $w$ is *irreducible*), and
3) There is no word $w_0 \in A^*$ spawning an infinite sequence of rewritings, $w_0 \to w_1 \to w_2 \to \cdots$. Such a system is called **Noetherian**.

We will say that $G$ **has a finite complete rewriting system** if there is a generating set $A$ for which there is a finite complete rewriting system for $(G, A)$.

We will need the following facts about finite complete rewriting systems.

**Proposition 0.** *The trivial group has a finite complete rewriting system.*

**Proposition 1.** $\mathbb{Z}$ *has a finite complete rewriting system.*

**Proposition 2 [Hr],[LC].** *If $G$ is a surface group, then $G$ has a finite complete rewriting system.*

**Proposition 3 [GS2].** *If $H$ is finite index in $G$ and $H$ has a finite complete rewriting system, then $G$ has a finite complete rewriting system.*

**Proposition 4 [GS1].** *If $1 \to K \to G \to Q \to 1$ is a short exact sequence, and $K$ and $Q$ have finite complete rewriting systems, then $G$ has a finite complete rewriting system.*

A thorough account of Thurston's eight geometries is given in [**Sc**]. If $M$ is a Riemannian manifold, then the Riemannian metric on $M$ lifts to a Riemannian metric on the universal cover, $\widetilde{M}$. Suppose now that $M$ is a closed 3-manifold with a uniform Riemannian metric. (This means that the isometry group of $\widetilde{M}$ acts transitively.) Thurston has shown that up to scaling, there are only eight possibilities for the Riemannian manifold $\widetilde{M}$ and that $\pi_1(M)$ is constrained in the following way. (We say that $G$ is **virtually** $H$ if $G$ contains a finite index copy of $H$.)

**Proposition 5.** *Suppose $M$ is a closed Riemannian 3-manifold with a uniform metric. Then one of the following holds:*

1) $\widetilde{M}$ is the 3-sphere and $\pi_1(M)$ is finite, i.e., virtually trivial.
2) $\widetilde{M}$ is Euclidean 3-space and $\pi_1(M)$ is virtually $\mathbb{Z}^3$.
3) $\widetilde{M}$ is $\mathbb{S}^2 \times \mathbb{R}$ and $\pi_1(M)$ is virtually $\mathbb{Z}$.
4) $\widetilde{M}$ is $\mathbb{H}^2 \times \mathbb{R}$ and $\pi_1(M)$ is virtually $H \times \mathbb{Z}$, where $H$ is the fundamental group of a closed hyperbolic surface.
5) $\widetilde{M}$ is Nil, the Lie group consisting of upper triangular real $3 \times 3$ matrices with one's on the diagonal, and $\pi_1(M)$ contains a finite index subgroup $G$ which sits in the short exact sequence $1 \to \mathbb{Z} \to G \to \mathbb{Z}^2 \to 1$.
6) $\widetilde{M}$ is $\widetilde{\mathrm{PSL}_2(\mathbb{R})}$, the universal cover of the unit tangent bundle of the hyperbolic plane, and $\pi_1(M)$ contains a finite index subgroup $G$ which sits in the short exact sequence $1 \to \mathbb{Z} \to G \to H \to 1$, where $H$ is the fundamental group of a closed hyperbolic surface.
7) $\widetilde{M}$ is Sol, a Lie group which is a semi-direct product of $\mathbb{R}^2$ with $\mathbb{R}$, and $\pi_1(M)$ has a finite index group $G$ which sits in the short exact sequence $1 \to \mathbb{Z}^2 \to G \to \mathbb{Z} \to 1$.
8) $\widetilde{M}$ is hyperbolic space. Under the further assumption that $M$ virtually fibers over a circle, $\pi_1(M)$ has a finite index group $G$ which sits in the short exact sequence $1 \to H \to G \to \mathbb{Z} \to 1$, where $H$ is the fundamental group of a closed hyperbolic surface.

The proof of Theorem 1 now consists of applying Propositions 1 – 4 to the cases of Proposition 5.

## §3. Non-uniform geometric 3-manifolds.

For an arbitrary closed 3-manifold $M$ satisfying Thurston's geometrization conjecture (see [**Sc**] for details), but not necessarily admitting a uniform Riemannian metric, finding rewriting systems becomes much more complicated. If $M$ is not orientable, then $M$ has an orientable double cover; Proposition 3 then says that if the fundamental group of the cover has a finite complete rewriting system, then so does $M$. So we may assume $M$ is orientable.

Any closed orientable 3-manifold $M$ can be decomposed as a connected sum $M = M_1 \# M_2 \# \cdots \# M_n$ in which each $M_i$ is either a closed irreducible 3-manifold, or is homeomorphic to $S^2 \times S^1$ ([**He**]). The fundamental group of $M$, then, can be written as the free product $\pi_1(M) = \pi_1(M_1) * \pi_1(M_2) * \cdots \pi_1(M_n)$. Another result of [**GS1**] says that the class of groups with finite complete rewriting systems is closed under free products; therefore, if $\pi_1(M_i)$ has a finite complete rewriting system for each $i$, then so does $\pi_1(M)$. So we may assume that our closed orientable 3-manifold is also irreducible.

Results of [**JS**] and [**Jo**] state that a closed irreducible 3-manifold $M$ can also be decomposed in a canonical way. There is a finite graph $\Gamma$ associated to $M$. For each vertex of $v$ of $\Gamma$, there is a compact 3-manifold $M_v \subset M$. For each edge $e$ of $\Gamma$ there is an incompressible torus $T_e^2 \subset M$. The boundary of $M_v$ consists of $\coprod_{v \in \partial e} T_e^2$ and $M$ is the union along these boundary tori of the pieces $M_v$. Consequently, the fundamental group of $M$ can also be realized as the group of the graph of groups given by placing the fundamental groups of the vertex manifolds at the corresponding vertices of $\Gamma$, and the fundamental group of a torus on each edge, together with the appropriate injections.

If $M$ satisfies Thurston's geometrization conjecture, then the interior of each one of these vertex manifolds admits a uniform Riemannian metric.

The simplest case of this type of decomposition occurs when the graph $\Gamma$ consists of a single vertex with no edges; this is dealt with in Proposition 5. In the case where $\Gamma$ has edges, the vertex manifolds are either cusped hyperbolic 3-manifolds or Seifert fibered manifolds with boundaries.

We take up the case in which $M$ is a graph of circle bundles. More specifically, suppose that $\Gamma$ is a finite graph, and suppose that at any vertex $v$ the vertex manifold $M_v$ is a circle bundle over a punctured surface with genus $g_v$. Then $M_v$ will have a torus boundary component for each of the punctures in the base surface; the edges of the graph $\Gamma$ determine how these boundary components will be glued together.

Recall that a **loop** is an edge with the same initial and terminal vertex. We need to assume that if all of the loops in the graph $\Gamma$ are removed, the resulting graph is a tree. Then the vertices of $\Gamma$ can be colored alternately red and blue, so that each edge that is not a loop joins a blue vertex to a red vertex. Orient all of the non-loop edges in $\Gamma$ by taking the blue vertex to be the initial vertex, so the red vertex is the terminal vertex. Let $V$ be the set of vertices of $\Gamma$, let $E$ be the set of edges in $\Gamma$ that are not loops, and let $L$ be the set of loops in $\Gamma$. For each edge $e \in E$, $\iota(e)$ will denote the initial vertex of $e$, and $\tau(e)$ will denote the terminal vertex; similarly for loops in $L$.

The fundamental group of the circle bundle at the vertex $v$ will be

$$\pi_1(M_v) = \langle a_{v1}, ..., a_{vg_v}, b_{v1}, ..., b_{vg_v}, \{p_e \mid e \in E, \iota(e) = v\}, \{q_e \mid e \in E, \tau(e) = v\},$$

$$\{r_l, s_l \mid l \in L, \iota(l) = v\} \mid \prod_{\iota(e)=v} p_e \prod_{\tau(e)=v} q_e \prod_{\iota(l)=v} r_l s_l = \prod_{j=1}^{g_v} [a_{vj}, b_{vj}]\rangle \times \langle x_v \rangle.$$

If $v$ is a blue vertex, then, there will be no generators of the form $q_e$, and if $v$ is a red vertex, there will be no generators of the form $p_e$ in the presentation for $\pi_1(M_v)$.

Suppose the edge $e \in E$ has initial vertex $v$ (so $v$ is blue) and terminal vertex $w$ (so $w$ is red). Then the amalgamation along this edge gives relations

$$x_v = q_e^{k_e} x_w^{n_e}$$
$$p_e = q_e^{k'_e} x_w^{n'_e}$$

where the matrix

$$\phi = \begin{pmatrix} k_e & n_e \\ k'_e & n'_e \end{pmatrix} \in \mathrm{SL}_2(\mathbb{Z}).$$

We are able to find finite complete rewriting systems in the case where this matrix is

$$\phi = \begin{pmatrix} 1 & n_e \\ 0 & 1 \end{pmatrix}.$$

Our relations in this case are

$$x_v = q_e x_w^{n_e}$$
$$p_e = x_w.$$

For each red vertex $w$, define the number $n_w$ to be the sum over all the edges $e \in E$, with target $\tau(e) = w$, of the numbers $n_e$.

If the loop $l$ has initial and terminal vertex $v$, then a generator $t_l$ is added, along with relations

$$t_l x_v t_l^{-1} = s_l^{k_l} x_v^{m_l}$$
$$t_l r_l t_l^{-1} = s_l^{k'_l} x_v^{m'_l}$$

where, again, the matrix

$$\phi = \begin{pmatrix} k_l & m_l \\ k'_l & m'_l \end{pmatrix} \in \mathrm{SL}_2(\mathbb{Z}).$$

As before, in order to find a finite complete rewriting system, we need to assume this matrix is

$$\phi = \begin{pmatrix} 1 & m_l \\ 0 & 1 \end{pmatrix};$$

our relations in this case are

$$t_l x_v t_l^{-1} = s_l x_v^{m_l}$$
$$t_l r_l t_l^{-1} = x_v$$

Replace the generator $p_e$ with the generator $x_{\tau(e)}$ in the presentation above, and replace the generator $q_e$ with the word $x_{\iota(e)} x_{\tau(e)}^{-n_e}$. The following is a rewriting system for the graph of circle bundles described above, with alphabet $A = S \cup S^{-1}$, where

$$S = \{x_v, a_{vj}, b_{vj}, r_l, s_l, t_l \mid v \in V, 1 \leq j \leq g_v, l \in L\}.$$

- inverse cancellation relators: $\{zz^{-1} \to 1, \quad z^{-1}z \to 1 \mid z \in S\}$
- blue vertex relators:
  $$\{x_v^{\pm 1} a_{vi}^{\pm 1} \to a_{vi}^{\pm 1} x_v^{\pm 1}, \qquad x_v^{\pm 1} b_{vi}^{\pm 1} \to b_{vi}^{\pm 1} x_v^{\pm 1},$$
  $$x_{\iota(k)}^{\pm 1} r_k^{\pm 1} \to r_k^{\pm 1} x_{\iota(k)}^{\pm 1}, \qquad x_{\iota(k)}^{\pm 1} s_k^{\pm 1} \to s_k^{\pm 1} x_{\iota(k)}^{\pm 1}\}$$
- red vertex relators:
  $$\{a_{wi}^{\pm 1} x_w^{\pm 1} \to x_w^{\pm 1} a_{wi}^{\pm 1}, \qquad b_{wi}^{\pm 1} x_w^{\pm 1} \to x_w^{\pm 1} b_{wi}^{\pm 1},$$
  $$r_l^{\pm 1} x_{\iota(l)}^{\pm 1} \to x_{\iota(l)}^{\pm 1} r_l^{\pm 1}, \qquad s_l^{\pm 1} x_{\iota(l)}^{\pm 1} \to x_{\iota(l)}^{\pm 1} s_l^{\pm 1}\}$$

- edge relators: $\{x_{\iota(e)}^{\pm 1} x_{\tau(e)}^{\pm 1} \to x_{\tau(e)}^{\pm 1} x_{\iota(e)}^{\pm 1}\}$
- blue amalgam relators:
$$\{a_{v1} b_{v1} \to \Lambda_v \prod_{j=g_v}^{2} [b_{vj}, a_{vj}] b_{v1} a_{v1},$$
$$a_{v1} b_{v1}^{-1} \to b_{v1}^{-1} \prod_{j=2}^{g_v} [a_{vj}, b_{vj}] \Lambda_v^{-1} a_{v1},$$
$$a_{v1}^{-1} \Lambda_v \prod_{j=g_v}^{2} [b_{vj}, a_{vj}] b_{v1} \to b_{v1} a_{v1}^{-1},$$
$$a_{v1}^{-1} b_{v1}^{-1} \to b_{v1}^{-1} a_{v1}^{-1} \Lambda_v \prod_{j=g_v}^{2} [b_{vj}, a_{vj}]\}$$
- red amalgam relators:
$$\{a_{w1} b_{w1} \to x_w^{-n_w} \Omega_w \prod_{j=g_w}^{2} [b_{wj}, a_{wj}] b_{w1} a_{w1},$$
$$a_{w1} b_{w1}^{-1} \to x_w^{n_w} b_{w1}^{-1} \prod_{j=2}^{g_w} [a_{wj}, b_{wj}] \Omega_w^{-1} a_{w1},$$
$$a_{w1}^{-1} \Omega_w \prod_{j=g_w}^{2} [b_{wj}, a_{wj}] b_{w1} \to x_w^{n_w} b_{w1} a_{w1}^{-1},$$
$$a_{w1}^{-1} b_{w1}^{-1} \to x_w^{-n_w} b_{w1}^{-1} a_{w1}^{-1} \Omega_w \prod_{j=g_w}^{2} [b_{wj}, a_{wj}]\}$$
- blue HNN relators:
$$\{x_{\iota(k)} t_k \to t_k r_k, \qquad x_{\iota(k)}^{-1} t_k \to t_k r_k^{-1},$$
$$r_k t_k^{-1} \to t_k^{-1} x_{\iota(k)}, \qquad r_k^{-1} t_k^{-1} \to t_k^{-1} x_{\iota(k)}^{-1},$$
$$s_k t_k \to t_k r_k^{-m_k} x_{\iota(k)}, \qquad s_k^{-1} t_k \to t_k r_k^{m_k} x_{\iota(k)}^{-1},$$
$$x_{\iota(k)} t_k^{-1} \to t_k^{-1} s_k x_{\iota(k)}^{m_k}, \qquad x_{\iota(k)}^{-1} t_k^{-1} \to t_k^{-1} s_k^{-1} x_{\iota(k)}^{-m_k}\},$$
and
- red HNN relators:
$$\{t_l r_l \to x_{\iota(l)} t_l, \qquad t_l r_l^{-1} \to x_{\iota(l)}^{-1} t_l,$$
$$t_l^{-1} x_{\iota(l)} \to r_l t_l^{-1}, \qquad t_l^{-1} x_{\iota(l)}^{-1} \to r_l^{-1} t_l^{-1},$$
$$t_l x_{\iota(l)} \to x_{\iota(l)}^{m_l} s_l t_l, \qquad t_l x_{\iota(l)}^{-1} \to x_{\iota(l)}^{-m_l} s_l^{-1} t_l,$$
$$t_l^{-1} s_l \to x_{\iota(l)} r_l^{-m_l} t_l^{-1}, \qquad t_l^{-1} s_l^{-1} \to x_{\iota(l)}^{-1} r_l^{m_l} t_l^{-1}\}$$

These rules range over all blue vertices $v$, red vertices $w$, and edges $e \in E$, as well as all loops $k$ at blue vertices and all loops $l$ at red vertices. The letter $\Lambda_v$ denotes the string of letters

$$\Lambda_v = \prod_{\iota(e)=v} x_{\tau(e)} \prod_{\iota(k)=v} r_k s_k.$$

$\Lambda_v^{-1}$, then, denotes the formal inverse of $\Lambda_v$, taking the letters in the string $\Lambda_v$ in the opposite order with their signs changed. The letter $\Omega_w$ denotes the string of letters

$$\Omega_w = \prod_{\tau(e)=w} x_{\iota(e)} \prod_{\iota(l)=w} r_l s_l,$$

and $\Omega_w^{-1}$ is its formal inverse.

Denote this set of rules to be $R$; the generators $A$ together with our rewriting rules $R$ give a presentation for the fundamental group of the graph of circle bundles.

**Theorem 2.** *The rewriting system $R$ on the set $A$ is a finite complete rewriting system for the fundamental group of the graph of circle bundles described above.*

*Proof.* In order to show that this rewriting system is complete, we will first show that a subset of the rules give rise to a complete rewriting system. Let

$$A' = A - \{t_k^{\pm 1} \mid k \in L, \iota(k) \text{ is blue}\},$$

and define $R'$ to be the rewriting system consisting of all of the rules above except the blue HNN relators and the inverse cancellation relators involving the letters of $A - A'$.

In order to show that this system $R'$ is Noetherian, we will show that there is a well-founded ordering on the words in $A'^*$ so that whenever a word is rewritten, the resulting word is smaller with respect to this order. This ordering is a recursive path ordering

**Definition [De].** *Let $>$ be a partial well-founded ordering on a set $A'$. The **recursive path ordering** $>_{rpo}$ on $A'^*$ is defined recursively from the ordering on $A'$ as follows. Given $s_1, ..., s_m, t_1, ..., t_n \in A'$, $s_1...s_m >_{rpo} t_1...t_n$ if and only if one of the following holds.*
*1) $s_1 = t_1$ and $s_2...s_m >_{rpo} t_2...t_n$.*
*2) $s_1 > t_1$ and $s_1...s_m >_{rpo} t_2...t_n$.*
*3) $s_2...s_m \geq_{rpo} t_1...t_n$.*
*The recursion is started from the ordering $>$ on $A'$ and from $s >_{rpo} 1$ for all $s \in A'$, where $1$ is the empty word in $A'^*$.*

Recursive path ordering is a well-founded ordering which is compatible with concatenation of words [**De**]. The following lemma is proved by inspection of the rules in the set $R'$.

**Lemma.** *Let $>$ be the recursive path ordering induced by*

$$t_l^{-1} > t_l > a_{w1}^{-1} > a_{w1} > b_{w1}^{-1} > b_{w1} > a_{w2}^{-1} > \cdots > b_{wg_w} > x_v^{-1} > r_l^{-1} > r_l > s_l^{-1} > s_l >$$
$$x_v > a_{v1}^{-1} > a_{v1} > b_{v1}^{-1} > b_{v1} > a_{v2}^{-1} > \cdots > b_{vg_v} > r_k^{-1} > r_k > s_k^{-1} > s_k > x_w^{-1} > x_w,$$

*where $v$ is any blue vertex, $w$ is any red vertex, $k$ is any loop at a blue vertex, and $l$ is any loop at a red vertex. Then for each of the rules $u \to v$ in $R'$, we have $u > v$.*

It follows from the Lemma that this system $R'$ is Noetherian. In order to show that $R'$ is also complete, it suffices to show that in the monoid presented by $(A', R')$, for each element $m$ in this monoid, there is exactly one word in $A'^*$ representing $m$ that cannot be rewritten.

The Knuth-Bendix algorithm [**KB**] is a computational procedure for checking that a Noetherian rewriting system is complete. This algorithm checks for overlapping for overlapping pairs of rules either of the form $r_1 r_2 \to s, r_2 r_3 \to t \in R'$ with $r_2 \neq 1$, or of the form $r_1 r_2 r_3 \to s, r_2 \to t \in R'$, where each $r_i \in A'^*$; these are called *critical pairs*. In the first case, the word $r_1 r_2 r_3$ rewrites to both $sr_3$ and $r_1 t$; in the second, it rewrites to both $s$ and $r_1 t r_3$. If there is a word $z \in A'^*$ so that $sr_3$ and $r_1 t$ both rewrite to $z$ in a finite number of steps in the first case, or so that $s$ and $r_1 t r_3$ both rewrite to $z$ in the second case, then the critical pair is said to be *resolved*. The Knuth-Bendix algorithm checks that all of the critical pairs of the system are resolved; if this is the case, then the rewriting system is complete. We have used this procedure to check that the rewriting system $R'$ is complete.

Since the rewriting system $R'$ is complete, for each word $u \in A'^*$, there is a bound on the lengths of all sequences of rewritings $u \to w_1 \to \cdots \to w_n$ (where the length of this sequence is defined to be $n$). The maximum of the lengths of all of the possible rewritings of $u$ is called the *disorder* of $u$, denoted $d_{R'}(u)$. We will use these numbers in order to show that the larger rewriting system $R$ is Noetherian.

In order to define a well-founded ordering on the set $A^*$, note that every word $w \in A^*$ can be written uniquely in the form

$$w = u_1 t_1 u_2 t_2 \cdots u_j t_j u_{j+1},$$

where each $u_i$ is a (possibly empty) word in $A'^*$ and each $t_i$ is a letter in $A - A'$. Define functions $\psi_i$ from $A^*$ to the nonnegative integers by

$$\begin{aligned}\psi_0(w) &= j, \\ \psi_{2i}(w) &= d_{R'}(u_i), \text{ and} \\ \psi_{2i+1}(w) &= length(u_i),\end{aligned}$$

where $i$ ranges from 1 to $j+1$, and *length* denotes the word length over $A'$. In order to compare words of different length, define $\psi_i(w) = 0$ if $i > 2j + 3$. For two words $w_1$ and $w_2$ in $A^*$, define $w_1 > w_2$ if $\psi_0(w_1) > \psi_0(w_2)$ or if $\psi_i(w_1) = \psi_i(w_2)$ for all $i < k$ and $\psi_k(w_1) > \psi_k(w_2)$. We claim this defines a well-founded ordering on $A^*$.

To check the claim, suppose $w \in A^*$. If a rule in $R'$ is applied to $w$, the rule must to be applied to one of the subwords $u_i$, so the value of $\psi_{2i}$ is reduced without altering the values of $\psi_k$ for any $0 \le k \le 2i - 1$. Suppose an inverse cancellation relator involving the letters of $A - A'$ is applied to $w$; in this case, the value of $\psi_0$ is reduced. Finally, if a blue HNN relator is applied to $w$, the rule must be applied to a subword $u_i t_i$ of $w$. Then the values of $\psi_k$ for any $0 \le k \le 2i - 1$ are not altered; the value of $\psi_{2i}$ either decreases or remains unchanged; and the value of $\psi_{2i+1}$ is reduced. So each time a word is rewritten, the resulting word is smaller with respect to this ordering. Therefore the rewriting system $R$ is also Noetherian.

Since $R$ is Noetherian, we have again applied the Knuth-Bendix procedure to check that the rewriting system $R$ is complete. □

## §4. An example

Rather than give the details of the Knuth-Bendix computation, we will give a description of the normal forms that these rewriting rules produce in the case where $M$ decomposes into two circle bundles. In this case, $G$ is a free product with amalgamation and the graph $\Gamma$ consists of two vertices $v$ and $w$ joined by a single edge. We assume that $v$ is blue and $w$ is red. $M_v$ and $M_w$ are circle bundles, each with a single torus boundary component and $M$ is formed by gluing along these tori. Thus we have

$$A = \pi_1(M_v) = \langle a_1, ..., a_g, b_1, ..., b_g \rangle \times \langle x \rangle$$
$$C = \pi_1(M_w) = \langle c_1, ..., c_h, d_1, ..., d_h \rangle \times \langle y \rangle$$
$$X = \langle x \rangle \times \langle y \rangle$$

so $G = \pi_1(M) = A *_X C$, where the gluing along $X$ is given by

$$x = \prod_{i=1}^{h} [c_i, d_i] y^n$$

and

$$\prod_{i=1}^{g} [a_i, b_i] = y.$$

The generating set for the fundamental group $G$ of $M$ will be $S \cup S^{-1}$ where

$$S = \{a_1, ..., a_g, b_1, ..., b_g, x, c_1, ..., c_h, d_1, ..., d_h, y\}.$$

The rewriting rules are:
- inverse cancellation relators: $\{zz^{-1} \to 1 \quad z^{-1}z \to 1 \mid z \in S\}$
- blue vertex relators: $\{x^{\pm 1} a_i^{\pm 1} \to a_i^{\pm 1} x^{\pm 1}, \quad x^{\pm 1} b_i^{\pm 1} \to b_i^{\pm 1} x^{\pm 1}\}$,
- red vertex relators: $\{c_i^{\pm 1} y^{\pm 1} \to y^{\pm 1} c_i^{\pm 1}, \quad d_i^{\pm 1} y^{\pm 1} \to y^{\pm 1} d_i^{\pm 1}\}$,
- edge relators: $\{x^{\pm 1} y^{\pm 1} \to y^{\pm 1} x^{\pm 1}\}$
- blue amalgam relators:
  $\{a_1 b_1 \to y p^{-1} b_1 a_1, \quad a_1 b_1^{-1} \to b_1^{-1} p y^{-1} a_1,$
  $a_1^{-1} y p^{-1} b_1 \to b_1 a_1^{-1}, \quad a_1^{-1} b_1^{-1} \to b_1^{-1} a_1^{-1} y p^{-1}\}$,
- red amalgam relators:
  $\{c_1 d_1 \to y^{-n} x q^{-1} d_1 c_1, \quad c_1 d_1^{-1} \to y^n d_1^{-1} q x^{-1} c_1,$
  $c_1^{-1} x q^{-1} d_1 \to y^n d_1 c_1^{-1}, \quad c_1^{-1} d_1^{-1} \to y^n d_1^{-1} c_1^{-1} x q^{-1}\}$

Here we use the letter $p$ to denote the string of letters

$$a_2 b_2 a_2^{-1} b_2^{-1} ... a_g b_g a_g^{-1} b_g^{-1},$$

so $p^{-1}$ denotes

$$b_g a_g b_g^{-1} a_g^{-1} ... b_2 a_2 b_2^{-1} a_2^{-1}.$$

Similarly the letter $q$ denotes the string

$$c_2 d_2 c_2^{-1} d_2^{-1} ... c_h d_h c_h^{-1} d_h^{-1},$$

and $q^{-1}$ denotes
$$d_h c_h d_h^{-1} c_h^{-1} \ldots d_2 c_2 d_2^{-1} c_2^{-1}.$$

To understand the normal forms that these produce, we consider several sublanguages.

Let $L(A/X)$ be the set of irreducible words on $\{a_i, b_i, y\}^{\pm 1}$ which do not end in $y^{\pm 1}$. Similarly, we let $L(X\backslash C)$ be the set of irreducible words on $\{c_i, d_i, x\}^{\pm 1}$ which do not begin in $x^{\pm 1}$. We take $L(X) = \{y^m x^n \mid m, n \in \mathbb{Z}\}$.

**Lemma.**
1) $L(A/X)$ bijects to $A/X$.
2) $L(X\backslash C)$ bijects to $X\backslash C$.
3) $L(X)$ bijects to $X$.

*Proof.* Clearly $L(A/X)$ surjects to $Z/X$, for the set of reduced words on these letters surjects to $A$, and deleting any trailing $y^{\pm 1}$ does not change the coset. Thus, to prove *1)* we must show that if $u, u' \in L(A/X)$ with $\overline{u}X = \overline{u'}X$ then $u = u'$. Here $u$ and $u'$ both evaluate into the free group on $\{a_i, b_i\}$, so in this case we have $\overline{uy^m} = \overline{u'}$ for some $m$. Observe that no rewriting rule has a left hand side consisting of letters of $\{a_i, b_i, y\}^{\pm 1}$ and ending in $y^{\pm 1}$ (other than free reduction). Thus, if $u$ is irreducible, then so is $uy^m$ for any $m$. If $m \neq 0$ we have two distinct irreducible words representing the same group element. However, it is not hard to carry out the Knuth-Bendix procedure on the set of rules evaluating into $A$. This ensures that for any element of $A$ there is a unique irreducible word and thus $u$ and $u'$ are identical as required.

The proof of *2)* is similar. Once again it is easy to see that $L(X\backslash C)$ surjects to $X\backslash C$. Now we suppose that $X\overline{w} = X\overline{w'}$ with $w$ and $w'$ in $L(X\backslash C)$. We then have $\overline{y^m x^n w} = \overline{w'}$ for some $m$ and $n$. Once again, these are both irreducible as there is no rewriting rule beginning in $y^m x^n$ that can be applied. Appealing to the Knuth-Bendix procedure in $C$ forces $m = n = 0$ and thus $w = w'$ as required.

The proof of *3)* is immediate. □

Now observe that any irreducible word $\theta$ has the form
$$\theta = u_1 v_1 w_1 \ldots u_k v_k w_k$$
where
- For each $i$, $v_i = y^{m_i} x^{n_i} \in L(X)$.
- For each $i$, $u_i$ is a maximal subword lying in $L(A/X)$.
- For each $i$, $w_i$ is a maximal subword lying in $L(X\backslash C)$.

The maximality of each $u_i$ ensures that no $v_i$ consists solely of $y^{\pm 1}$'s directly preceding $v_{i+1}$. Likewise the maximality of each $w_i$ ensures that no $v_i$ consists solely of $x^{\pm 1}$'s directly following $v_{i-1}$.

We call $k$ the **length** of $\theta$. Let $L_k$ be the set of all irreducible words of length $k$. For each $g \in G$ the *AC-length* of $g$ is the minimal $k$ such that $g \in (AC)^k$. Let $G_k$ be the set of all group elements of $AC$-length $k$. That is, $g \in G_k$ if $k$ is minimal such that $g = A_1 C_1 \ldots A_k C_k$ with $A_i \in A$, $C_i \in C$.

**Claim.** $L_k$ *bijects to* $G_k$.

*Proof.* This is an induction on $k$.

We check the case $k = 1$. Clearly $L_1$ surjects to $AC$, since any element of $A$ has the form $u_1 v'_1$ and each element of $C$ has the form $v''_1 w_1$. Multiplying these together and applying the replacement rules produces a word of the form $u_1 v_1 w_1$ as required.

We now check that the map from $L_1$ to $AC$ is injective. Suppose $g \in AC$ and $g = \overline{u_1 v_1 w_1}$ Notice that $A/X$ bijects to $AC/C$. Thus $g$ determines a coset $gC$ in $AC/C$ and thus a unique element of $A/X$. Consequently, $g$ determines $u_1$.

On the other hand, having determined the coset representative $u_1$ of $gC$ in $AC/C$, there is a unique $c \in C$ so that $g = \overline{u_1} c$ and this, in turn, determines $v_1 w_1$.

We now assume by induction that $L_k$ bijects to $G_k$ and check that $L_{k+1}$ bijects to $G_{k+1}$. It is easy to see that $L_{k+1}$ surjects to $G_{k+1}$. For suppose $g \in G_{k+1}$. Then $g$ has the form $g_k h$ with $g_k \in G_k$, $h \in G_1$. We represent $g_k$ by a word of $L_k$ and $h$ by a word of $L_1$. We concatenate these words and apply our rewriting rules. The resulting word $\theta$ lies in $\cup_{i=1}^{k+1} L_i$. Since $\cup_{i=1}^{k} L_i$ misses $G_{k+1}$, it follows that $\theta \in L_{k+1}$.

We must show that $L_{k+1}$ injects to $G_{k+1}$. Notice that $g_k$ and $h$ are determined by $g$ up to an element of $X$. Thus, if $g = g_k h = g'_k h'$ then there is $z \in X$ so that $g'_k = g_k z$ and $h' = z^{-1} h$. Suppose then that $g$ is represented by two irreducible words

$$\theta = u_1 v_1 w_1 \ldots u_k v_k w_k u_{k+1} v_{k+1} w_{k+1}$$
$$\theta' = u'_1 v'_1 w'_1 \ldots u'_k v'_k w'_k u'_{k+1} v'_{k+1} w'_{k+1}$$

We take $g_k$ and $g'_k$ to be the group elements represented by the $L_k$ portions of $\theta$ and $\theta'$. Thus $h$ and $h'$ are represented by the remaining portions of these two words.

Notice that if $w_k$ ends in $x^{\pm 1}$ then $u_{k+1}$ and $v_{k+1}$ are both empty, for otherwise, the final $x^{\pm 1}$'s of $w_k$ would have moved right through $u_{k+1}$ and any $y^{\pm 1}$'s of $v_{k+1}$. This cannot happen, since each $w_i$ was chosen to be maximal. In the same manner, we do not have $w_k$ empty and $v_k$ ending in $x^{\pm 1}$. The same argument applies to $\theta'$. Suppose the element $z \in X$ is represented by the word $y^m x^n \in L(X)$. Then if the word $u_1 v_1 w_1 \ldots u_k v_k w_k y^m x^n$ is reduced using the rules of the rewriting system, for the resulting irreducible word $u_1 v_1 w_1 \ldots u_{k-1} v_{k-1} w_{k-1} \tilde{u}_k \tilde{v}_k \tilde{w}_k$, we have that either $\tilde{w}_k = w_k x^n$ or $\tilde{w}_k = w_k$ is empty and $\tilde{v}_k = v_k x^n$. Now this irreducible word represents the same element of $G_k$ as $u'_1 v'_1 w'_1 \ldots u'_k v'_k w'_k$. Therefore our induction hypothesis says that $\tilde{v}_k = v'_k$ and $\tilde{w}_k = w'_k$. So either $\tilde{w}_k = w_k x^n = w'_k$ or else $\tilde{w}_k$ and $w'_k$ are both empty and $\tilde{v}_k = v_k x^n = v'_k$. Since $w'_k$ cannot end with $x^{\pm 1}$ and if $w'_k$ is empty then $v'_k$ cannot end with $x^{\pm 1}$, this shows that $n$ must be zero. Consequently, $u_1 v_1 w_1 \ldots u_k v_k w_k$ and $u'_1 v'_1 w'_1 \ldots u'_k v'_k w'_k$ differ by at most a power of $y$ in $G$; that is, we have $z = \overline{y^m}$.

On the other hand, if $u_{k+1}$ begins in $y^{\pm 1}$, we must have $v_k$ and $w_k$ empty, for any leading $y^{\pm 1}$'s of $w_{k+1}$ would have had to move left through $w_k$ and any $x^{\pm 1}$'s of $v_k$. Maximality of the words $u_i$ does not allow this to happen. Also, we cannot have $u_{k+1}$ is empty and $v_{k+1}$ beginning with $y^{\pm 1}$. It follows by a similar argument, then, that $u_{k+1}v_{k+1}w_{k+1}$ and $u'_{k+1}v'_{k+1}w'_{k+1}$ differ in $G$ by at most a power of $x$; that is, $z = \overline{x^n}$. Since $z$ is now both a power of $x$ and a power of $y$, that power is plainly 0, so $g_k = g'_k$ and $h = h'$. By induction

$$u_1 v_1 w_1 \ldots u_k v_k w_k = u'_1 v'_1 w'_1 \ldots u'_k v'_k w'_k$$

and

$$u_{k+1} v_{k+1} w_{k+1} = u'_{k+1} v'_{k+1} w'_{k+1},$$

so $\theta = \theta'$ as required. □

Since $G = \coprod G_k$ and the language of irreducible words is $\coprod L_k$ it follows that the language of irreducible words is a normal form which bijects to $G$.

## §5. A question

When the gluings of the circle bundles at the vertices of $\Gamma$ are more complicated, or when the circle bundles themselves are replaced by more general Seifert-fibered spaces, we were unable to find finite complete rewriting systems. So we end with the following.

**Question.** *Does every fundamental group of a closed 3-manifold satisfying Thurston's geometrization conjecture have a finite complete rewriting system?*

## Acknowledgement

In the course of our work we have used Rewrite Rule Laboratory [**KZ**], a software package for performing the Knuth-Bendix algorithm, to check completeness of our rewriting systems for many examples.

New Mexico State University,
Las Cruces, NM 88003
USA

University of Melbourne
Parkville, Vic 3052
Australia